# О неускоренных эффективных методах решения разреженных задач квадратичной оптимизации

*Аникин А.С. (ИДСТУ СО РАН)*

*Гасников А.В. (ИППИ РАН; ПреМоЛаб МФТИ)*

*Горнов А.Ю. (ИДСТУ СО РАН)*

**Аннотация**

В работе исследуется класс разреженных задач квадратичной оптимизации и его обобщения, на которых будут эффективно работать специальные модификации обычного (неускоренного) прямого градиентного метода с l1-нормой и метода условного градиента (Франк–Вульфа).

**Ключевые слова:** huge-scale оптимизация, разреженность, l1-норма, метод Франк–Вульфа.

## 1. Введение

В работе [1] (см. также [2]) были предложены два различных способа решения задачи минимизации (на единичном симплексе) квадратичного функционала специального вида, связанного с задачей PageRank [3]. Первый способ базировался на обычном (неускоренном) прямом градиентном методе в l1-норме [4]. Второй способ базировался на методе условного градиента [5] (Франк–Вульф). Общей особенностью обоих способов является возможность учитывать разреженность задачи в большей степени, чем известные альтернативные способы решения отмеченной задачи. За счет такого свойства, итоговые оценки времени работы упомянутых методов в ряде случаев получаются лучше, чем оценки для ускоренных методов и их различных вариантов [2, 6]. Естественно было задаться вопросами: на какие классы задач и на какие простые множества (входящие в ограничения) можно перенести разработанные в [1] алгоритмы? Настоящая статья имеет целью частично ответить на поставленные вопросы.

В п. 2 мы показываем (во многом следуя изначальным идеям Ю.Е. Нестерова [7, 8]), что прямой градиентный метод в l1-норме из [1] можно распространить на общие задачи квадратичной минимизации. При этом можно перейти от задач на симплексе к задачам на всем пространстве (это даже заметно упростит сам метод и доказательство оценок). Но наиболее интересно то, что удается распространить этот метод на общие аффинно-сепарабельные задачи с сепарабельными композитами [6]. Такие задачи часто возникают при моделировании компьютерных и транспортных сетей [9–11], а также в анализе данных [12]. Основными конкурирующими методами здесь являются различные варианты (прямые, двойственные) ускоренных покомпонентных методов [6, 7, 9, 13, 14], которые чаще интерпретируются как различные варианты методов рандомизации суммы [6, 12, 15–22]. Таким образом, в п. 2 предлагается новый подход к этим задачам.

В п. 3 мы распространяем метод Франк–Вульфа из [1] на более общий класс задач квадратичной минимизации. К сожалению, выйти за пределы описанного в п. 3 класса задач нам не удалось. Тем не менее, стоит отметить интересную (саму по себе) конструкцию, позволившую перенести метод с единичного симплекса на неотрицательный ортант. Стоит также отметить, что эта конструкция существенным образом опиралась на ряд новых идей, успешно использовавшихся и в других контекстах.

В п. 4 мы предлагаем новый рандомизированный метод (дополнительно предполагая сильную выпуклость функционала задачи), который работает быстрее неускоренного покомпонентного метода, и позволяет для ряда задач в большей степени воспользоваться их разреженностью, чем позволяют это сделать ускоренные покомпонентные методы.

"За кадром" описываемых далее сюжетов стоит задача разработки эффективных численных методов решения системы линейных уравнений $Ax = b$ в пространстве огромной размерности (от десятков миллионов переменных и выше). При этом качество решения оценивается согласно $\|Ax - b\|_2$. Полученные в работе результаты позволяют надеяться, что предложенные в пп. 2–4 методы будут доминировать остальные в случае, когда (приводим наиболее важные условия):

1) Матрица $A$ такова, что в каждом столбце и каждой строке не более $s \ll \sqrt{n}$ отличных от нуля элементов.

2) Решение $x^*$ системы $Ax = b$ таково, что $\|x^*\|_1 \approx \|x^*\|_2$ ну или точнее: всегда имеет место неравенство

$$\|x^*\|_2 \le \|x^*\|_1 \le \sqrt{n}\|x^*\|_2,$$

которое (в нашем случае) можно уточнить следующим образом

$$\|x^*\|_2 \le \|x^*\|_1 \ll \sqrt{n}\|x^*\|_2.$$

Сравнительному анализу предложенных в данной работе методов с альтернативными методами поиска решения системы $Ax = b$ (в частности с популярными сейчас ускоренными покомпонентными методами [6]) посвящено приложение.

## 2. Прямой неускоренный градиентный метод в l1-норме

Рассматривается следующая задача

$$f(x) = \frac{1}{2}\langle Ax, x \rangle - \langle b, x \rangle \to \min_{x \in \mathbb{R}^n}, \qquad (1)$$

где квадратная матрица $A$ предполагается симметричной неотрицательно определенной. Также предполагается, что в каждом столбце и каждой строке матрицы $A$ не более $s \ll n$ элементов отлично от нуля. Мы хотим решить задачу (1) с точностью $\varepsilon$

$$f(x^N) - f(x^*) \le \varepsilon.$$

Заметим, что тогда

$$\|Ax^N - b\|_2 \le \sqrt{2\lambda_{\max}(A)\varepsilon}.$$

Эту задачу предлагается решать обычным градиентным методом, но не в евклидовой норме, а в l1-норме [4] или (что то же самое для данной задачи) неускоренным покомпонентным методом с выбором максимальной компоненты [7] (если $i^k$ определяется не единственным образом, то можно выбрать любого представителя):

$$x^{k+1} = x^k + \arg\min_{h \in \mathbb{R}^n}\left\{f(x^k) + \langle \nabla f(x^k), h \rangle + \frac{L}{2}\|h\|_1^2\right\} = \begin{cases} x_i^k, & i \ne i^k = \arg\max_{j=1,\ldots,n}\left|\partial f(x^k)/\partial x_j\right| \\ x_i^k - \frac{1}{L}\frac{\partial f(x^k)}{\partial x_{i^k}}, & i = i^k \end{cases}, \qquad (2)$$

где $L = \max_{i,j=1,\ldots,n}|A_{ij}|$. Точку старта $x^0$ итерационного процесса выберем таким образом, чтобы только одна из компонент была отлична от нуля. Любая итерация такого метода (за исключением самой первой – первая требует $\mathrm{O}(n)$ арифметических операций) может быть осуществлена за $\mathrm{O}(s\log_2 n)$. Действительно, $x^{k+1}$ и $x^k$ отличаются в одной компоненте. Следовательно,

$$\nabla f(x^{k+1}) = Ax^{k+1} - b = Ax^k - b + \frac{1}{L}\frac{\partial f(x^k)}{\partial x_{i^k}}Ae_{i^k} = \nabla f(x^k) + \frac{1}{L}\frac{\partial f(x^k)}{\partial x_{i^k}}Ae_{i^k},$$

где

$$e_i = \underbrace{(0,...,0,\overset{n}{1},0,...,0)}_{i}.$$

Таким образом, по условию задачи (матрица $A$ разрежена) $\nabla f(x^{k+1})$ отличается от $\nabla f(x^k)$ не более чем в $s$ компонентах. Следовательно, $i^{k+1}$ можно пересчитать за $O(s\log_2 n)$, храня массив $\left\{\left|\partial f(x^k)/\partial x_j\right|\right\}_{j=1}^n$ в виде специальной структуры данных (кучи) для поддержания максимального элемента [1, 8].

Согласно [4] необходимое число итераций

$$N = \frac{2\max\limits_{i,j=1,...,n}|A_{ij}|R_1^2}{\varepsilon},$$

где

$$R_1^2 = \sup\left\{\|x-x^*\|_1^2:\ f(x) \le f(x^0)\right\}.$$

Таким образом, общее число арифметических операций (время работы метода) можно оценить следующим образом

$$O\left(n + s\log_2 n \frac{\max\limits_{i,j=1,...,n}|A_{ij}|R_1^2}{\varepsilon}\right),$$

Рассмотрим теперь задачу вида

$$f(x) = \sum_{k=1}^m f_k(A_k^T x) + \sum_{k=1}^n g_k(x_k) \to \min_{x \in \mathbb{R}^n},$$

где все функции (скалярного аргумента) $f_k$, $g_k$ имеют равномерно ограниченные числом $L$ первые две производные,[1] а матрица $A = [A_1,...,A_m]^T$ такова, что в каждом ее столбце не больше $s_m \ll m$ ненулевых элементов, а в каждой строке – не больше $s_n \ll n$. Описанный выше подход (формула (2)) позволяет решить задачу за время

$$O\left(n + s_n s_m \log_2 n \frac{L\max\limits_{i=1,...,n}\|A^{\langle i\rangle}\|_2^2 R_1^2}{\varepsilon}\right),$$

---

[1] $g_k''$ можно равномерно ограничивать числом $L\max\limits_{i=1,...,n}\|A^{\langle i\rangle}\|_2^2$.

где $A^{\langle i \rangle}$ – $i$-й столбец матрицы $A$.

Собственно, задача из [1] является частным случаем приведенной выше общей постановки (в смысле выбора функционала):

$$f(x) = \frac{1}{2}\|Ax\|_2^2 + \frac{\gamma}{2}\sum_{i=1}^{n}(-x_i)_+^2 \to \min_{\langle x,e \rangle = 1}.$$

В принципе, можно перенести результаты, полученные в этом пункте на случай, когда оптимизация происходит на симплексе. Это частично (но не полностью [1]) решает основную проблему данного метода: большое значение $R_1^2$, даже когда $\|x^*\|_1^2$ небольшое. Однако мы не будем здесь этого делать.

## 3. Метод условного градиента (Франк–Вульфа)

Рассмотрим теперь следующую задачу[2]

$$f(x) = \frac{1}{2}\langle Ax, x \rangle - \langle b, x \rangle \to \min_{x \in \mathbb{R}_+^n}, \tag{3}$$

где квадратная матрица $A$ предполагается симметричной неотрицательно определенной. Также предполагаем, что в каждом столбце и каждой строке матрицы $A$ не более $s \ll n$ элементов отлично от нуля, и в векторе $b$ не более $s$ элементов отлично от нуля.

Предположим сначала, что мы знаем такой $R$, что решение задачи (3) удовлетворяет условию

$$x^* \in \overline{S}_n(R) = \left\{ x \in \mathbb{R}_+^n : \sum_{i=1}^{n} x_i^* \leq R \right\}.$$

Выберем одну из вершин $\overline{S}_n(R)$ и возьмем точку старта $x^1$ в этой вершине. Далее будем действовать по индукции, шаг которой имеет следующий вид.

Решаем задачу

$$\langle \nabla f(x^k), y \rangle \to \min_{y \in \overline{S}_n(R)}. \tag{4}$$

---

[2] Отличие от задачи (1) в том, что рассматривается неотрицательный ортант вместо всего пространства – раздутием исходного пространства в два раза к такому ограничению можно прийти из задачи оптимизации на всем пространстве.

Введем (если $i_k$ определяется не единственным образом, то можно выбрать любого представителя)

$$i_k = \arg\min_{i=1,\ldots,n} \partial f(x^k)/\partial x_i.$$

Обозначим решение задачи (4) через

$$y^k = \begin{cases} R \cdot e_{i_k}, & \text{если } \partial f(x^k)/\partial x_{i_k} < 0 \\ 0, & \text{если } \partial f(x^k)/\partial x_{i_k} \geq 0 \end{cases}.$$

Положим

$$x^{k+1} = (1-\gamma_k)x^k + \gamma_k y^k, \ \gamma_k = \frac{2}{k+1}, \ k=1,2,\ldots,$$

Имеет место следующая оценка [5, 23, 24]

$$f(x^N) - f(x^*) \leq f(x^N) - \max_{k=1,\ldots,N}\left\{f(x^k) + \langle \nabla f(x^k), y^k - x^k\rangle\right\} \leq \frac{2L_p R_p^2}{N+1}, \qquad (5)$$

где

$$R_p^2 = \max_{x,y \in \overline{S}_n(R)} \|y-x\|_p^2, \ L_p = \max_{\|h\|_p \leq 1} \langle h, Ah\rangle, \ 1 \leq p \leq \infty.$$

С учетом того, что оптимизация происходит на $\overline{S}_n(R)$, мы выбираем $p=1$. Несложно показать, что этот выбор оптимален. В результате получим, что

$$R_1^2 = 4R^2, \ L_1 = \max_{i,j=1,\ldots,n}|A_{ij}|. \qquad (6)$$

Таким образом, чтобы

$$f(x^N) - f(x^*) \leq \varepsilon,$$

достаточно сделать

$$N = \frac{8 \max_{i,j=1,\ldots,n}|A_{ij}| R^2}{\varepsilon}$$

итераций.

Сделав один раз дополнительные вычисления стоимостью $\mathrm{O}(n)$, можно так организовать процедуру пересчета $\nabla f(x^k)$ и вычисление $i_k$, что каждая итерация будет иметь сложность $\mathrm{O}(s\log_2 n)$. Для этого вводим

$$\beta_k = \prod_{r=1}^{k-1}(1-\gamma_r), \ z^k = x^k/\beta_k, \ \tilde{\gamma}_k = \gamma_k/\beta_{k+1}.$$

Тогда итерационный процесс можно переписать следующим образом

$$z^{k+1} = z^k + \tilde{\gamma}_k y^k.$$

Пересчитывать $Az^{k+1} - b$ при известном значении $Az^k$ можно за $\mathrm{O}(s)$. Далее, задачу

$$i_{k+1} = \arg\min_{i=1,\ldots,n} \partial f(x^{k+1})/\partial x_i$$

можно переписать как

$$i_{k+1} = \arg\min_{i=1,\ldots,n}\left(\left[Az^{k+1}\right]_i - \frac{b_i}{\beta_{k+1}}\right).$$

Поиск $i_{k+1}$ можно осуществить за $\mathrm{O}(s\log_2 n)$ (см. п. 2). Определив в конечном итоге $z^N$, мы можем определить $x^N$, затратив дополнительно не более

$$\mathrm{O}(n + \ln N) = \mathrm{O}(n)$$

арифметических операций. Таким образом, итоговая оценка сложности описанного метода будет

$$\mathrm{O}\left(n + s\log_2 n \frac{\max_{i,j=1,\ldots,n}|A_{ij}|R_1^2}{\varepsilon}\right),$$

где $R_1^2 = \|x^*\|_1^2$ (следует сравнить с аналогичной формулой из п. 2).

Отметим, что при этом мы можем пересчитывать (это не надо делать, если известно значение $f(x^*)$, например, для задачи (3) естественно рассматривать постановки с $f(x^*) = 0$)

$$f(x^{k+1}) + \langle \nabla f(x^{k+1}), y^{k+1} - x^{k+1}\rangle$$

также за $\mathrm{O}(s\log_2 n)$. Следуя [23], можно немного упростить приведенные рассуждения за счет небольшого увеличения числа итераций. А именно, можно использовать оценки (при этом следует полагать $\gamma^k = 2(k+2)^{-1}$)

$$f(x^k) - f(x^*) \leq \langle \nabla f(x^k), x^k - y^k \rangle,$$

$$\min_{k=1,\ldots,N} \langle \nabla f(x^k), x^k - y^k \rangle \leq \frac{7L_1 R_1^2}{N+2}.$$

Вернемся теперь к предположению, что изначально известно $R$. На практике, как правило, даже если мы можем как-то оценить $R$, то оценка получается слишком завышенной. Используя формулы (5), (6) мы можем воспользоваться следующей процедурой рестартов по параметру $R$.

Выберем сначала $R = R^0 = 1$. Делаем предписанное этому $R$ число итераций и проверяем (без дополнительных затрат) критерий останова (все необходимые вычисления уже были сделаны по ходу итерационного процесса)

$$f\left(x^N\right) - \max_{k=1,\ldots,N}\left\{f\left(x^k\right) + \left\langle\nabla f\left(x^k\right), y^k - x^k\right\rangle\right\} \le \frac{8\max_{i,j=1,\ldots,n}|A_{ij}|R^2}{N+1} \le \varepsilon.$$

Если он выполняется, то мы угадали и получили решение. Если не выполняется, то полагаем $R := \chi R$ ($\chi > 1$), и повторяем рассуждения. Остановимся поподробнее на вопросе оптимального выбора параметра $\chi$ [25]. Обозначим $R^* = \|x^*\|_1$. Пусть

$$R^0 \chi^{r-1} < R^* \le R^0 \chi^r.$$

Общее число итераций будет пропорционально

$$1 + \chi^2 + \chi^4 + \chi^{2r} = \frac{\chi^{2r+2}-1}{\chi^2-1} \le \frac{\chi^4}{\chi^2-1}\left(\frac{R^*}{R^0}\right)^2.$$

Выберем $\chi = \sqrt{2}$, исходя из минимизации правой части этого неравенства. При этом общее число итераций возрастет не более чем в четыре раза по сравнению со случаем, когда значение $R^*$ известно заранее.

Описанный выше подход распространяется и на задачи

$$f(x) = \frac{1}{2}\|Ax - b\|_2^2 \to \min_{x \in \mathbb{R}_+^n},$$

Матрица $A$ такова, что в каждом ее столбце не больше $s_m \ll m$ ненулевых элементов, а в каждой строке – не больше $s_n \ll n$. В векторе $b$ не более $s_m$ элементов отлично от нуля. Описанный выше подход (на базе метода Франк–Вульфа) позволяет решить задачу за время

$$O\left(n + s_n s_m \log_2 n \frac{\max_{i=1,\ldots,n}\left\|A^{\langle i\rangle}\right\|_2^2 R_1^2}{\varepsilon}\right),$$

где $R_1^2 = \|x^*\|_1^2$ (следует сравнить с аналогичной формулой из п. 2, полученной для более общего класса задач).[3]

В связи написанным выше в этом пункте, заметим, что задача может быть не разрежена, но свойство разреженности появляется в решении при использовании метода Франк–Вульфа, что также может заметно сокращать объем вычислений. Интересные примеры имеются в работах п. 3.3 [26], [27].

В последнее время методы условного градиента переживают бурное развитие в связи с многочисленными приложениями к задачам машинного обучения. В связи с этим появились интересные обобщения классического метода Франк–Вульфа. Отметим, например, работы [28, 29]. Интересно было бы понять: возможно ли перенести (а если возможно, то как именно и в какой степени) результаты этого пункта на более общий класс задач (чем класс задач с квадратичным функционалом) и на более общий класс методов?

---

[3] Здесь также как и в задаче (3) требуются рестарты по неизвестному параметру $R_1$, который явно используется в методе в качестве размера симплекса, фигурирующего при решении вспомогательной задачи ЛП на каждой итерации. Однако, как было продемонстрировано выше, все это приводит к увеличению числа итераций не более чем в четыре раза.

## 4. Неускоренный рандомизированный градиентный спуск в сильно выпуклом случае

Рассмотрим для большей наглядности снова задачу (1). Будем считать, что $f(x)$ $\mu$-сильно выпуклый функционал в 2-норме, т.е. $\lambda_{\min}(A) \geq \mu > 0$. Сочетая методы [30] и [31], введем следующий рандомизированный метод

$$x^{k+1} = x^k - \frac{2}{\mu \cdot (k+1)} \|\nabla f(x^k)\|_1 \operatorname{sign}\left(\frac{\partial f(x^k)}{\partial x_{i(x^k)}}\right) e_{i(x^k)}, \ x^1 = 0,$$

где

$$e_i = (\underbrace{0,...,0,1,0,...,0}_{i})^{\overbrace{\phantom{0,...,0,1,0,...,0}}^{n}},$$

$$i(x^k) = i \text{ с вероятностью } \frac{1}{\|\nabla f(x^k)\|_1} \left|\frac{\partial f(x^k)}{\partial x_i}\right|, \ i = 1,...,n.$$

Тогда[4] [30]

$$E[f(y^k)] - f(x^*) \leq \frac{\frac{2}{N}\sum_{k=1}^{N} E\left[\|\nabla f(x^k)\|_1^2\right]}{\mu \cdot (N+1)} \leq \frac{2 \max_{k=1,...,N} E\left[\|\nabla f(x^k)\|_1^2\right]}{\mu \cdot (N+1)},$$

$$E\left[\|\nabla f(y^k)\|_2^2\right] \leq \frac{4L \max_{k=1,...,N} E\left[\|\nabla f(x^k)\|_1^2\right]}{\mu \cdot (N+1)},$$

где $L = \lambda_{\max}(A)$, а

$$y^N = \frac{2}{N \cdot (N+1)} \sum_{k=1}^{N} k \cdot x^k.$$

---

[4] Второе неравенство может быть достаточно грубым [6].

На базе описанного метода построим новый метод, который "следит" за последовательностью $\sum_{k=1}^{N} k \cdot x^k$, пересчитывая (с некоторой частотой)

$$\left\| \nabla f\left( \frac{2}{N \cdot (N+1)} \sum_{k=1}^{N} k \cdot x^k \right) \right\|_1^2.$$

Метод "дожидается" момента[5] $N = \mathrm{O}(nL/\mu)$, когда

$$\left\| \nabla f\left( y^N \right) \right\|_1^2 \leq \frac{1}{2} \left\| \nabla f\left( x^1 \right) \right\|_1^2.$$

и перезапускается с $x^1 := y^N$. Можно показать, что необходимое число таких перезапусков для достижения точности $\varepsilon$ (по функции) будет $\sim \log_2(\varepsilon^{-1})$. Более того, подобно [6], можно так организовать описанную выше процедуру, чтобы попутно получить и оценки вероятностей больших отклонений (детали мы вынуждены здесь опустить).

Предположим теперь, что $\nabla f(x + he_i)$ отличается от $\nabla f(x)$ (при произвольных $x$, $h$ и $e_i$) не более чем в $s$ компонентах ($s \ll n$).[6] Для задачи (1) это имеет место (можно также говорить и о задаче (3), с очевидной модификацией описанного метода – все оценки сохраняются). Можно так организовать процедуру выбора момента $N$ (с сохранением свойства $N = \mathrm{O}(nL/\mu)$), что амортизационная (средняя) сложность итерация метода (с учетом затрат на проверку условий остановки на каждом перезапуске) будет $\mathrm{O}(s \log_2 n)$ [1, 3]. Казалось бы, что мы ничего не выиграли по сравнению с обычными (неускоренными) покомпонентными методами [6]. Количество итераций и стоимость одной итерации одинаковы для обоих методов (в категориях O( ) с точностью до логарифмических множителей). В действительности, ожидается, что предложенный нами алгоритм будет работать заметно быстрее (неускоренного) покомпонентного метода, поскольку, как уже отмечалось, при получении этих оценок мы пару раз использовали

---

[5] Это довольно грубая оценка на $N$, поскольку использующееся при ее получении правое неравенство

$$1 \leq \left\| \nabla f(y^N) \right\|_1^2 \Big/ \left\| \nabla f(y^N) \right\|_2^2 \leq n$$

может быть сильно завышенным.

[6] Это условие можно обобщить с сохранением оценок, например, на случай, когда $f(x) = g(x) + h(x)$, и существует такая (скалярная) функция $\alpha(x, he_i) > 0$, что $\nabla f(x + he_i) = \nabla g(x + he_i) + \nabla h(x + he_i)$ отличается от $\alpha(x, he_i) \cdot \nabla g(x) + \nabla h(x)$ не более чем в $s$ компонентах. Для этого приходится перезаписать исходный метод: отличие в том, что теперь вводится двойная рандомизация (рандомизация согласно вектору $\nabla g(x)$ и независимая рандомизация согласно вектору $\nabla h(x)$).

потенциально довольно грубые неравенства. С другой стороны пока мы говорили только о задаче (1) (в разреженном случае), для которой эффективно работают ускоренные покомпонентные методы [6] с такой же оценкой стоимости одной итерации, но заметно лучшей оценкой для числа итераций $N = \mathrm{O}\left(n\sqrt{L/\mu}\right)$ (с точностью до логарифмического множителя). В связи с этим может показаться, что предложенный в этом пункте метод теоретически полностью доминируем ускоренными покомпонентными методами. На самом деле, это не совсем так. Хорошо известно, см., например, [6], что существующие сейчас всевозможные модификации ускоренных покомпонетных методов, которые могут сполна учитывать разреженность задачи (что проявляется в оценке стоимости итерации $\mathrm{O}(s)$), применимы лишь к специальному классу задач [6]. Для общих задач стоимость одной итерации ускоренного покомпонентного метода будет $\mathrm{O}(n)$ независимо от разреженности задачи. Например, на данный момент не известны такие модификации ускоренных покомпонентных методов, которые бы позволяли сполна учитывать разреженность в задаче [6, 25, 31]

$$f(x) = \ln\left(\sum_{i=1}^{m}\exp\left(A_k^T x\right)\right) - \langle x,b\rangle + \frac{\mu}{2}\|x\|_2^2 \to \min_{x\in\mathbb{R}^n}.$$

В то время как описанный выше метод распространим и на эту задачу с оценкой амортизационной стоимости одной итерации $\mathrm{O}(s_m s_n \log_2 n)$.

### Приложение

В этом приложении мы напомним некоторые полезные факты о поиске решения (псевдорешения) системы линейных уравнений $Ax = b$ [32, 33]. К близким задачам и методам приводит изучение проекции точки на аффинное множество $Ax = b$ [9, 31].

Хорошо известно, что задача поиска такого вектора $x$, что $Ax = b$, может быть эффективно полиномиально решена даже в концепции битовой сложности (например, простейшим методом Гаусса или более современными методами внутренней точки). Однако требуется, учитывая специфику задачи, так подобрать метод, чтобы можно было искать решения систем огромных размеров. Несмотря на более чем двухвековую историю, эта область математики (эффективные численные методы решения систем линейных уравнений) до сих пор бурно развивается. О чем, например, говорит недавний доклад Д. Спильмана на международном математическом конгрессе [35].

Пусть наблюдается вектор

$$b = Ax + \varepsilon,$$

где матрица $A$ – известна, $\varepsilon_k \in N(0, \sigma^2)$ – независимые одинаково распределенные случайные величины, $k = 1,...,m$ (ненаблюдаемые). Оптимальная оценка[7] неизвестного вектора параметров $x$ определяется решением задачи

$$f(x) = \frac{1}{2}\|Ax - b\|_2^2 \to \min_{x \in \mathbb{R}^n}. \qquad (7)$$

Введем псевдообратную матрицу (Мура–Пенроуза) $A^\dagger = (A^T A)^{-1} A^T$, если столбцы матрицы $A$ линейно независимы (первый случай) и $A^\dagger = A^T (AA^T)^{-1}$, если строки матрицы $A$ линейно независимы (второй случай). В первом случае $x^* = A^\dagger b$ – единственное решение задачи (7), во втором случае $x^* = A^\dagger b$ – решение задачи (7) с наименьшим значением 2-нормы.[8] Вектор $x^*$ называют псевдорешением задачи $Ax = b$.

Поиск решения системы $Ax = b$ был сведен выше к решению задачи (7). Для задачи выпуклой оптимизации (7) существуют различные эффективные численные методы. Например, метод сопряженных градиентов, сходящийся со скоростью (стоимость одной итерации этого метода равна $\mathrm{O}(sn)$ – числу ненулевых элементов в матрице $A$)

$$\frac{1}{2}\|Ax^N - b\|_2^2 = f(x^N) - f(x^*) \le \varepsilon,$$

$$N = \mathrm{O}\left(\sqrt{\frac{\lambda_{\max}(A^T A)\|x^*\|_2^2}{\varepsilon}}\right),$$

где $N \le n$. На практике при не очень больших размерах матрицы $A$ ($\max\{m, n\} \le 10^6$) эффективно работают квазиньютоновские методы типа L-BFGS [34]. А при совсем небольших размерах ($\max\{m, n\} \le 10^3 - 10^4$) и методы внутренней точки [33]. Однако нам интересны ситуации, когда $\max\{m, n\} \gg 10^6$. В таких случаях эффективно использовать прямые и двойственные ускоренные покомпонентные методы для задачи (7). Эти методы дают следующие оценки [6] (одна итерация у обоих методов в среднем требует $\mathrm{O}(s)$ арифметических операций)[9]

$$E\left[\frac{1}{2}\left\|Ax^N - b\right\|_2^2\right] \le \varepsilon, \qquad \text{(прямой метод)}$$

$$N = \mathrm{O}\left(n\sqrt{\frac{\bar{L}_x R_x^2}{\varepsilon}}\right), \; \bar{L}_x^{1/2} = \frac{1}{n}\sum_{k=1}^{n}\left\|A^{\langle k \rangle}\right\|_2, \; R_x = \left\|x^*\right\|_2,$$

$$E\left[\left\|Ax^N - b\right\|_2\right] \le \varepsilon, \qquad \text{(двойственный метод)}$$

$$N = \mathrm{O}\left(n\sqrt{\frac{\bar{L}_y R_y}{\varepsilon}}\right), \; \bar{L}_y^{1/2} = \frac{1}{m}\sum_{k=1}^{m}\left\|A_k\right\|_2, \; R_y = \left\|y^*\right\|_2,$$

где $y^*$ – решение (если решение не единственно, то можно считать, что $y^*$ – решение с наименьшей 2-нормой) "двойственной" задачи

$$-\frac{1}{2}\left\|A^T y\right\|_2^2 + \langle b, y \rangle \to \max_{y \in \mathbb{R}^m}.$$

---

[7] То есть несмещенная и с равномерно (по $x$) минимальной дисперсией. При такой (статистической) интерпретации также стоит предполагать, что $n \le m$.

[8] В приложениях во втором случае чаще ищут не псевдорешение, а вводят в функционал регуляризатор [9], например, исходя из байесовского подхода (дополнительного априорного вероятностного предположения относительно вектора неизвестных параметров) или исходя из желания получить наиболее разреженное решение. В таком случае многое (но не все) из того, что написано в данной статье удается сохранить (должным образом модифицировав). При этом появляются и дополнительные новые возможности за счет появления новых степеней свободы в выборе регуляризатора [27].

Напомним, что $\left(A_k\right)^T$ – $k$-я строка матрицы $A$, а $A^{\langle k \rangle}$ – $k$-й столбец.

Из написанного выше не очевидно, что для получения решения $x^*$ необходимо исходить именно из решения задачи оптимизации (7) и, как следствие, исходить из определяемого этой задачей критерия точности решения (2-норма невязки). Например, в настоящей работе мы исходили в основном из задачи (1), которая хотя и похожа по свойствам на задачу (7), но все же приводит к отличному от (7) критерию точности решения. Вместо задачи (7) в определенных ситуациях вполне осмысленно бывает рассматривать даже негладкую (но выпуклую) задачу

$$f(x) = \left\|Ax - b\right\|_\infty \to \min_{x \in \mathbb{R}^n}.$$

Для такой задачи можно предложить специальным образом рандомизированный вариант метода зеркального спуска [36], который гарантируют выполнение (считаем $f(x^*) = 0$)

$$E\left[\left\|Ax^N - b\right\|_\infty\right] \le \varepsilon,$$

за время

---

[9] Далее для наглядности мы дополнительно предполагаем, что $f(x^*) = 0$. Если $f(x^*) = f^* > 0$, то из того что $f(x^N) - f^* \le \varepsilon$ будет следовать $\left\|Ax^N - b\right\|_2 - f^* \le \varepsilon/f^*$ (этот простое, но полезное наблюдение было сделано Дмитрием Островским (INRIA Гренобль)).

$$\mathrm{O}\left( s_m n + s_m \log_2 m \frac{\max\limits_{k=1,\ldots,m} \|A_k\|_1^2 R_x^2}{\varepsilon^2} \right),$$

если в каждом столбце матрицы $A$ не больше $s_m \le m$ ненулевых элементов. Отметим [36], что обычный (не рандомизированный) метод зеркального спуска гарантирует выполнение

$$\|Ax^N - b\|_\infty \le \varepsilon$$

за время

$$\mathrm{O}\left( m + n + s_n s_m \log_2 m \frac{\max\limits_{k=1,\ldots,m} \|A_k\|_2^2 R_x^2}{\varepsilon^2} \right),$$

если в каждом столбце матрицы $A$ не больше $s_m \le m$ ненулевых элементов, а в каждой строке – не больше $s_n \le n$.

Выше мы обсуждали, так называемые, вариационные численные методы, которые сводили поиск решения системы $Ax = b$, к решению задачи выпуклой оптимизации. Однако имеется большой класс итерационных численных методов (типа метода простой итерации), которые исходят из перезаписи системы $Ax = b$ в эквивалентном виде. Например, в виде $x = \tilde{A}x + b$, где $\tilde{A} = I - A$, и организации вычислений по формуле[10] $x^{k+1} = \tilde{A}x^k + b$ [39]. Для ряда таких методов также удается получить оценки скорости сходимости, причем не только для нормы невязки в системе, но и непосредственно для оценки невязки в самом решении $\|x^N - x^*\|_2$. К сожалению, такие оценки получаются весьма пессимистичными, причем эти оценки типично реализуются на практике (закон Мерфи). Поясним сказанное примером, восходящим к работе Красносельского–Крейна [40]. Пусть

---

[10] Интересно заметить, что критерием сходимости такого итерационного процесса является: $\lambda_{\max}(\tilde{A}) < 1$ [32]. Однако, если при этом $\|\tilde{A}\|_{2,2} = \sigma_{\max}(\tilde{A}) = \lambda_{\max}(A^T A) > 1$ (к счастью, такие ситуации не типичны [37]), то существуют такие $x^0$ (причем на них довольно легко попасть), стартуя с которых, итерационная процедура на первых итерациях приводит к резкому росту нормы невязки $\|x^N - x^*\|_2$, т.е. к наличию "горба" у такой последовательности. Высота этого горба может быть (в зависимости от конкретной задачи) сколь угодно большой [32, 38]. Поскольку все вычисления проходят с конечной длиной мантиссы, то такой резкий рост может приводить к вычислительной неустойчивости процедуры, что неоднократно наблюдалось в численных экспериментах [32]. Если $\sigma_{\max}(\tilde{A}) < 1$, то таких проблем нет, и последовательность $\|x^N - x^*\|_2$ мажорируется геометрической прогрессией с основанием $\sigma_{\max}(\tilde{A})$.

$$x^{k+1} = \tilde{A}x^k + b,$$

где $\tilde{A} = I - A$ – положительно определенная матрица $n \times n$ с собственными числами $0 < \lambda_1 \leq ... \leq \lambda_n < 1$, а $x^0$ – выбирается равновероятно из $B_R(0)$ – 2-шара радиуса $R$ с центром в $0$. Обозначим через

$$\delta^k = x^{k+1} - x^k = \tilde{A}x^k + b - x^k = b - Ax^k$$

контролируемую невязку. В качестве критерия останова итерационного процесса выберем момент $N$ первого попадания $\delta^N \in B_\varepsilon(0)$. Заметим, что

$$\left\| x^k - x^* \right\|_2 = \left\| \left( I - \tilde{A} \right)^{-1} \delta^k \right\|_2 = \left\| A^{-1} \delta^k \right\|_2.$$

Имеет место следующий результат. При $R \to \infty$ вероятность выполнения неравенств

$$0.999 \cdot \frac{\varepsilon \lambda_n}{1 - \lambda_n} \leq \left\| x^N - x^* \right\|_2 \leq \frac{\varepsilon}{1 - \lambda_n}$$

стремится к единице. Поскольку типично, что при больших $n$ число $\lambda_n$ близко к 1, то вряд ли можно рассчитывать на то, что $\left\| x^N - x^* \right\|_2$ мало, если установлена только малость $\left\| Ax^N - b \right\|_2$.

Отметим, что оценки скорости сходимости многих итерационных методов (типа простой итерации) исходят из первого метода Ляпунова, т.е. из спектральных свойств матрицы $A$. При плохих спектральных свойствах вводится регуляризация. По-сути, большинство таких методов можно проинтерпретировать как численный метод решения соответствующей вариационной переформулировки задачи. Регуляризация при этом имеет четкий смысл (восходящий к пионерским работам А.Н. Тихонова) – сделать функционал сильно выпуклым и использовать его сильную выпуклость при выборе метода и получении оценки скорости сходимости [41]. Однако вариационный подход позволяет сполна использовать и второй метод Ляпунова, что дает возможность распространить всю мощь современных численных методов выпуклой оптимизации на решение задачи $Ax = b$. В частности, использовать рандомизацию [36] или(и) свойства неравнозначности компонент в решении $x^*$ и разреженность $A$, которые использовались при разработке методов из настоящей работы.

Ю.Е. Нестеровым недавно было отмечено, что возможные прорывы в разработке новых эффективных численных методов (в том числе применительно к решению системы $Ax = b$) можно ожидать от рандомизированных квазиньютоновских методов. К сожалению, здесь имеются большие проблемы, унаследованные от обычных (не рандомизированных) квазиньютоновских методов, с получением оценок скоростей сходимости, адекватных реальной (наблюдаемой на практике) скорости сходимости (все известные сейчас оценки никак не объясняют быстроту сходимости этих методов на практике). Тем не менее, недавно появилось несколько интересных работ Гаверса–Ричтарика в этом направлении [42–44].

Авторы выражают благодарность Ю.Е. Нестерову за внимательное отношение к работе.

## Литература